\documentclass[a4paper,9pt]{article}
\usepackage{amsfonts}
\usepackage{amssymb,amsthm,color,tikz}
\usepackage{amsmath}
\usepackage{comment}
\usepackage{dsfont}

\usepackage[normalem]{ulem}
\usepackage{url}

\newcommand{\R}{\mathbb{R}}
\newcommand{\N}{\mathbb{N}}

\newtheorem{theorem}{Theorem}

\newtheorem{corollary}[theorem]{Corollary}

\newtheorem{definition}{Definition}

\def\dd {\mathrm{d}}
\newcommand{\Om}{\Omega}

\newcommand{\bp}{\begin{proof}}
\newcommand{\ep}{\end{proof}}

\begin{document}
\title{Qualitative properties of the heat content}

\author{{Michiel van den Berg} \\
School of Mathematics, University of Bristol\\
Fry Building, Woodland Road\\
Bristol BS8 1UG, United Kingdom\\
\texttt{mamvdb@bristol.ac.uk}\\
\\
{Katie Gittins}\\
Department of Mathematical Sciences, Durham University\\
Mathematical Sciences \& Computer Science Building\\
Upper Mountjoy Campus, Stockton Road\\
Durham DH1 3LE, United Kingdom\\
\texttt{katie.gittins@durham.ac.uk}}
\date{22 April 2025}\maketitle
\vskip 2truecm \indent
\begin{abstract}\noindent
We obtain monotonicity and convexity results for the heat content of domains in Riemannian manifolds and in Euclidean space subject to various initial temperature conditions.
We introduce the notion of a \emph{strictly decreasing temperature set}, and show that it is a sufficient condition to ensure monotone heat content.
In addition, in Euclidean space, we construct a domain and an initial condition for which the heat content is not monotone, as well as a domain and an initial condition for which the heat content is monotone but not convex.
\end{abstract}
\vskip 1truecm \noindent   {Mathematics Subject
Classification ({2020})}: { 35K05, { 35K20 or 35K15}.\\ \textbf{Keywords}: Heat content, convexity, monotonicity.

\section{Introduction and main results \label{sec1}}
\label{sec1}

In this paper, we investigate qualitative properties of heat flow problems in open sets in Euclidean space and in Riemannian manifolds without boundary conditions.

For example, consider an open set $\Omega\subset\R^m$ that is initially at temperature $1$ while its complement, $\R^m \setminus \Omega$, is initially at temperature $0$.
No boundary conditions are imposed on the boundary $\partial \Om$ of $\Omega$ and the heat equation evolves on $\R^m$.

A version of the isoperimetric inequality for the heat semigroup corresponding to this heat flow problem has been established in \cite{P03} by making use of the connection between the perimeter of the set and the small-time asymptotic behaviour of the semigroup (see also \cite{MPPP}).

We study the interplay between the geometry of $\Om$ and the heat content of $\Om$, that is, the amount of heat left inside $\Om$ at time $t$.
The refined asymptotic behaviour of the heat content of $\Om$ as $t\downarrow 0$ has been obtained in a variety of geometric settings. For example, polygons in $\R^2$ \cite{MvdBKG2}, horn-shaped regions in $\R^m$ \cite{MvdB1}, and smooth, compact Riemannian manifolds contained in a larger compact Riemannian manifold \cite{MvdBPG1}.
Two-sided bounds for the heat content and for the heat loss were obtained in \cite{MvdBKG1} for the case of an open set in $\R^m$ with $R$-smooth boundary and finite Lebesgue measure, and in \cite{MvdB2} for the case of an open set in a complete, smooth, non-compact, $m$-dimensional Riemannian manifold.
More recently the heat content has been analysed in the context of metric measure spaces and sub-Riemannian manifolds \cite{C, Rossi}.

The goal of this paper is to investigate the monotonicity and convexity of the heat content of $\Om$ as a function of $t$, with various initial data, in the setting where no boundary conditions are imposed on $\partial \Om$.

%We introduce the following notation and definitions, and recall the discussion from p.885-886 in \cite{MvdBGR} to set the scene for our main results.

Let $M$ be a smooth, connected, complete and stochastically complete $m$-dimensional Riemannian manifold and let $\Delta $ be the Laplace-Beltrami
operator acting on functions in $L^2(M)$. It is well known (see \cite{EBD3}, \cite{GB})
that the heat equation
\begin{equation}\label{e1}
\Delta u=\frac{\partial u}{\partial t},\quad x\in M,\quad t>0,
\end{equation}
has a unique, minimal, positive fundamental solution $p_M(x,y;t)$
where $ x\in M$, $y \in M$, $t>0$. This solution, called the Dirichlet
heat kernel for $M$, is symmetric in $x,y$, strictly positive,
jointly smooth in $x,y\in M$ and $t>0$, and it satisfies the
semigroup property
\begin{equation}\label{e2}
p_M(x,y;s+t)=\int_{M} \dd z\,\ p_M(x,z;s)p_M(z,y;t) ,
\end{equation}
for all $x,y\in M$ and $t,s>0$, where $\dd z $ is the Riemannian
measure on $M$. In addition
\begin{equation}\label{e2b}
	\int_M \dd y\, p_M(x,y;t)= 1
\end{equation}
since $M$ is stochastically complete.
Let $\Omega$ be an open subset of $M$. Equation (\ref{e1}) with the initial condition
\begin{equation}\label{e3}
u(x;0^+)=\psi (x),\quad x\in \Omega,
\end{equation}
has a solution
\begin{equation}\label{e4}
u_{\Omega,\psi }(x;t)=\int_{\Om} \dd y\,\ p_M(x,y;t)\psi(y),
\end{equation}
for any function $\psi$ on $\Om$ from a variety of function spaces. For example, let $\psi\in C_b(\Om),\psi\ge 0$, $\psi \not\equiv 0$, the set of bounded continuous functions from $\Om$ into $[0,\infty).$
Then initial condition (\ref{e3}) is understood in the sense that $u_{\Om,\psi}\left(
\cdot;t\right) \rightarrow \psi \left( \cdot\right) $ as $t\downarrow 0$, where
the convergence is locally uniform.

\bigskip
Let $\Omega$ be a non-empty, open subset of $M$, and let $\psi:\Omega\rightarrow [0,\infty)$ be bounded and measurable.
We define the heat content of $\Omega$ with initial datum $\psi$ by
\begin{equation}\label{e5}
H_{\Omega,\psi}(t)=\int_{\Omega}\int_{\Omega}\dd x\,\dd y\, p_M(x,y;t)\psi(y).
\end{equation}

It was shown in \cite[Proposition 1]{P} that if $\Omega \subset \R^m$ is bounded, then $t \mapsto H_\Omega(t)$ is decreasing and convex. In Theorem \ref{the1} below we consider the more general situation of a Riemannian manifold.
A particular case of interest is when $\psi \equiv 1$ on $\Om$ for which we write
\begin{equation*}%\label{e6}
H_{\Om}(t)=H_{\Om,1}(t).
\end{equation*}

We introduce the following definition that will give us a sufficient condition that ensures monotone heat content.

\begin{definition}\label{def1}
	Let $M$ be a smooth, connected, complete and stochastically complete $m$-dimensional Riemannian manifold.
	An open set $\Om\subset M$ is a (strictly) decreasing temperature set if for all $x\in \Om$, $t\mapsto u_{\Om,1}(x;t)$ is (strictly) decreasing.
\end{definition}

Our first main result is the following.
\begin{theorem}\label{the1}
	Let $\Omega$ be a non-empty, open subset of $M$, where $M$ is a smooth, connected, complete and stochastically complete $m$-dimensional Riemannian manifold.
\begin{enumerate}
\item[\textup{(i)}]If $H_{\Om}(t)<\infty$ for all $t>0$, then $t\mapsto H_{\Om}(t)$ is decreasing and convex. Moreover, $\lim_{t \to \infty} H_\Om (t)$ exists.
\item[\textup{(ii)}]If $H_{\Om}(t)<\infty$ for all $t>0$ and if $\lim_{t \to \infty} H_\Om (t)=0$, then all right-hand derivatives of $H_{\Om}(t)$ with respect to $t$ are strictly negative, and $t\mapsto H_{\Om}(t)$ is strictly decreasing.
\item[\textup{(iii)}]If $M$ is in addition closed, then $t\mapsto H_{\Omega}(t)$ is strictly decreasing, and strictly convex if and only if $|M\setminus \Om|>0$.
\item[\textup{(iv)}]If $\Om\subset M$ is a (strictly) decreasing temperature set with finite measure, and if $\psi : \Om \to [0,\infty)$ is bounded and measurable, then $t\mapsto H_{\Om,\psi}(t)$ is (strictly) decreasing.
\end{enumerate}
\end{theorem}

All remaining theorems concern results for Euclidean space $\R^m$, for which
\begin{equation}\label{ee}
	p_{\R^m}(x,y;t)= \frac{e^{-\frac{|x-y|^2}{4t}}}{(4\pi t)^{m/2}}.
	\end{equation}

For the case where $\psi \equiv 1$ on $\Om \subset \R^m$, it was shown in \cite[Proposition 8]{MvdB1} that
if $\Om$ is convex, then $\Om$ is a decreasing temperature set. In Example 6 of \cite{MvdB1} it was shown that the disjoint union of a ball and a suitable concentric annulus in $\R^2$ is not a decreasing temperature set.
Below we show that the disjoint union of two balls with equal radii $\delta$ in $\R^m$ at distance $2$ is a strictly decreasing temperature set for some $\delta$ sufficiently small. So the convexity assumption in  \cite[Proposition 8]{MvdB1} is sufficient but not necessary and sufficient.

\begin{theorem}\label{the3}  Let $\Omega_{\delta}= B_{\delta}(c_1)\cup B_{\delta}(c_2)\subset \R^m$, $m \in \mathbb{N}$, and $c_1=(-1-\delta,0,...,0),\,c_2=(1+\delta,0,...,0)$. Let $\psi = {\bf 1}_{\Om_\delta}$.
	If $\delta=\frac{1}{20}$, then $\Omega_{\delta}$ is a strictly decreasing temperature set.
\end{theorem}

\medskip

In Theorem \ref{the5} below, for non-empty, open, bounded sets in $\R^m$, we obtain a lower bound for the second derivative of the heat content and show that this derivative is bounded away from $0$ for all $t$ sufficiently large, and an upper bound for the first derivative of the heat content and show that this derivative is bounded away from $0$.

\begin{theorem}\label{the5}If  $\Omega$ is a non-empty, open set in $\R^m$ with $\textup{diam}(\Om)<\infty$, then
	\begin{itemize}
		\item[\textup{(i)}]\begin{equation}\label{e49}
			\frac{d^2H_{\Om}(t)}{dt^2}\ge \frac{4m^2+4m-7}{16t^2}H_{\Om}(t),\,t\ge (\textup{diam}(\Om))^2,
		\end{equation}
		\item[\textup{(ii)}]\begin{equation}\label{e50}
			\frac{d^2H_{\Om}(t)}{dt^2}\ge \begin{cases}0,\,t < (\textup{diam}(\Om))^2,\\
				\\ (4m^2+4m-7)\pi^2e^{-1/4}\frac{|\Om|^2}{(4\pi t)^{(m+4)/2}},\,t\ge (\textup{diam}(\Om))^2,
			\end{cases}
		\end{equation}
		\item[\textup{(iii)}]\begin{equation}\label{e51}
			\frac{dH_{\Om}(t)}{dt}\le \begin{cases} -\frac{4m^2+4m-7}{2(m+2)}\pi e^{-1/4}\frac{|\Om|^2}{(4\pi (\textup{diam}(\Om))^2)^{(m+2)/2}}, 0<t\le (\textup{diam}(\Om))^2,\\
				\\
				-\frac{4m^2+4m-7}{2(m+2)}\pi e^{-1/4}\frac{|\Om|^2}{(4\pi t)^{(m+2)/2}}, t\ge (\textup{diam}(\Om))^2,
			\end{cases}
		\end{equation}
		\item[\textup{(iv)}]\begin{equation*}%\label{e52}
			\frac{dH_{\Om}(t)}{dt}\le -\frac{4m^2+4m-7}{8(m+2)t} e^{-1/4}H_{\Om}(t),\,t\ge (\textup{diam}(\Om))^2.
		\end{equation*}
	\end{itemize}
\end{theorem}

Theorem \ref{the5} can be generalised to the case of non-negative, measurable initial temperature $\psi$ as follows.
\begin{corollary}\label{cor1}
	If  $\Omega$ is a non-empty, open set in $\R^m$ with $\textup{diam}(\Om)<\infty$ and if $\psi\ge 0$, $\psi \not\equiv 0$, bounded and measurable, then analogous results to those of Theorem \ref{the5} hold for $H_{\Om, \psi}$. Moreover, the analogues of parts \textup{(i), (iv)} hold with $H_{\Om}$ replaced by $H_{\Om, \psi}$ and the analogues of parts \textup{(ii), (iii)} hold with $|\Om |^2$ replaced by $|\Om | \int_{\Om} \psi(y) \, dy$.
\end{corollary}
\noindent
Corollary \ref{cor1} follows immediately from the proof of Theorem \ref{the5}.

\medskip

We now explore the effects of changing the initial datum $\psi$ on the monotonicity and convexity of $t \mapsto H_{\Om, \psi}(t)$.

Throughout for $r_2>r_1>0,\,\tilde{c}\in \R^m$, we let $B_{r_1}{ (\tilde{c})}=\{x\in\R^m: |x-\tilde{c}|<r_1\}$, $A_{(r_1,r_2)}=\{x\in\R^m:r_1<|x|<r_2\}$, and $\omega_m=|B_1(0)|$.

First, we construct an example which shows that if $\psi$ is not constant on $\Omega$, then $t \mapsto H_{\Om,\psi}(t)$ need not be monotone in $t$.

\begin{theorem}\label{the2} Let $m \in \mathbb{N}, \, c>2,\,t^*>e^{9/(2m)}$, and let
\begin{equation*}%\label{e26}
\Omega_c=B_1(0) \cup A_{(2, c)} \subset \R^m,\,\,\,	\psi(x)={\bf 1}_{B_1{ (0)}}(x).
\end{equation*}
If
\begin{equation}\label{e25}
c>(32t^*)^{1/2}\Big(\log\Big(\frac{2^{5m/2} \Gamma((m+2)/2)}{2^m-1}\big(e^{-9/4}-{t^*}^{-m/2}\big)^{-1}\Big)\Big)^{1/2},	
\end{equation}
then $H_{\Om_c,\psi}(0)> H_{\Om_c,\psi}(t^*)>H_{\Om_c,\psi}(1),$ so that  $t\mapsto H_{\Om_c,\psi}(t)$ is not monotone.
\end{theorem}

As a consequence of Theorem \ref{the2}, since $t\mapsto H_{\Om_c,\psi}(t)$ is not monotone, we deduce by Theorem \ref{the1}(iv) that, for $c$ sufficiently large, $\Om_c$ is not a decreasing temperature set.

\bigskip
In addition, we construct an example which shows that if $\psi$ is not constant on $\Omega$, then $t \mapsto H_{\Om,\psi}(t)$ can be a decreasing function of $t$ but need not be convex.
\begin{theorem}\label{the4}
	Let $\Omega=B_1(0)\subset\R^m$, $m \in \mathbb{N}$, $\psi(y)=|1-|y||^{\alpha}$ with $\alpha>1$, then $t\rightarrow H_{B_1(0),\psi}(t)$ is decreasing but not convex.
\end{theorem}

\noindent The proofs of Theorems \ref{the1}, \ref{the3}, \ref{the5}, \ref{the2}, \ref{the4} are deferred to Sections \ref{sec2}, \ref{sec4}, \ref{sec6}, \ref{sec3}, \ref{sec5} respectively.

\section{Proof of Theorem \ref{the1} \label{sec2}}
\begin{proof}
\noindent(i) For $t>0$, by \eqref{e4} we have
\begin{equation}\label{e7}
u_{\Omega}(x;t)=\int_{\Omega}\dd y\,p_M(x,y;t).
\end{equation}

We first show that $t\mapsto H_{\Om}(t)$ is decreasing. By \eqref{e2} and \eqref{e7} we have for $t>0,s>0$,
\begin{align}\label{e8}
u_{\Omega}(x;t+s)&=\int_{\Omega}\dd y\,p_M(x,y;t+s)\nonumber\\&
=\int_{\Omega}\dd y\, \int_{M}\dd z\, p_M(x,z;t)p_M(z,y;s)\nonumber\\&
=\int_{M}\dd z\, p_M(x,z;t)u_{\Omega}(z;s),
\end{align}
where we have used Tonelli's Theorem in the last identity. Integrating \eqref{e8} with respect to $x$ over $\Omega$ yields
\begin{equation}\label{e9}
H_{\Omega}(t+s)=\int_{M}\dd x\,u_{\Omega}(x;t)u_{\Omega}(x;s).
\end{equation}
By \eqref{e9}, \eqref{e8}, \eqref{e2}, \eqref{e2b} and symmetry of the heat kernel
\begin{align*}%\label{e10}
H_{\Omega}(t+s)&=\int_{M}\dd x\,u_{\Omega}(x;(t+s)/2)^2\nonumber\\&
=\int_M\dd x\,\int_M\dd y_1\,p_M(x,y_1;s/2)u_{\Om}(y_1;t/2)\int_Mdy_2\,p_M(x,y_2;s/2)u_{\Om}(y_2;t/2)\nonumber\\&
=\int_M\dd y_1\,\int_M\dd y_2\,p_M(y_1,y_2;s)u_{\Om}(y_1;t/2)u_{\Om}(y_2;t/2)\nonumber\\&\le
\frac12\int_M\dd y_1\,\int_M\dd y_2\,p_M(y_1,y_2;s)\Big(u_{\Om}(y_1;t/2)^2+u_{\Om}(y_2;t/2)^2\Big)\nonumber\\&=
\int_Mdy_1\,\int_M\dd y_2\,p_M(y_1,y_2;s)u_{\Om}(y_1;t/2)^2\nonumber\\&
 = \int_M\dd y_1\,u_{\Om}(y_1;t/2)^2\nonumber\\&
=H_{\Omega}(t).
\end{align*}

\smallskip

To prove convexity, we first note that since $H_{\Om}(t)<\infty,\,t>0$, it suffices to prove that $H$ is midpoint convex. See pp.164--167 in \cite{RPB}.
Let $t>0,\, \delta>0$. By \eqref{e9}, we have
\begin{align*}%\label{e11}
\frac12\Big(H_{\Om}(t)+H_{\Om}(t+2\delta)\Big)&=\frac12\int_{M}\dd z\,\Big(u_{\Om}(z;t/2)^2+u_{\Om}(z;(t+2\delta)/2)^2\Big)\nonumber\\&
\ge \int_{M}\dd z\,u_{\Om}(z;t/2)u_{\Om}(z;(t+2\delta)/2)\nonumber\\&
=H_{\Om}(t+\delta).
\end{align*}
This proves the convexity of $H_{\Om}$. Since the map $t \mapsto H_\Om(t)$ is decreasing and bounded from below, $\lim_{t \to \infty} H_\Om(t)$ exists.

\medskip

\noindent(ii) By convexity of $t\mapsto H_{\Om}(t)$, we have that the right-hand derivative $H_{\Om}^{'+}(t):= \lim_{\varepsilon\downarrow 0}\varepsilon^{-1}({H_{\Om}(t+\varepsilon)-H_{\Om}(t)})$ is non-decreasing in $t$.
See p.167 in \cite{RPB}. Hence if $H_{\Om}^{'+}(T)\ge 0$ for some $T>0$, then $H_{\Om}^{'+}(t)\ge 0,\,t\ge T.$ This in turn implies $H_{\Om}(t)\ge H_{\Om}(T)>0,\,t\ge T.$ This contradicts $\lim_{t \to \infty} H_\Om (t)=0,$
and $H_{\Om}^{'+}(T)< 0,\,T>0$. Hence $t\mapsto H_{\Om}(t)$ is strictly decreasing.

\medskip

\noindent(iii) We note that if $|M\setminus\Om|=0$, then
\begin{align*}%\label{e16}
	H_{\Om}(t)= &\int_{\Om}\int_{\Om}dx\,dy\,p_{M}(x,y;t)\nonumber\\&
	=\int_{M}\int_{M}dx\,dy\,p_{M}(x,y;t)\nonumber\\&
	=H_M(t)=|M|,
\end{align*}
since smooth, closed Riemannian manifolds are stochastically complete. Hence $t\mapsto H_{\Omega}(t)$ is constant, and is not strictly decreasing nor is it strictly convex.
Next consider the case $0<|\Om|<|M|$. Since $M$ is smooth and closed, the spectrum of the Laplace-Beltrami operator $\Delta$ acting in $L^2(M)$ is discrete and consists of eigenvalues
$\{\mu_1(M)\le \mu_2(M)\le ...\}$ accumulating at $\infty$ only. Let $\{u_{j,M},\,j\in \N\}$ denote a corresponding orthonormal basis of eigenfunctions. Since $M$ is connected,
$\mu_1(M)=0$ and has multiplicity $1$. Furthermore $u_{1,M}=|M|^{-1/2}$. The minimal heat kernel for $M$ has an $L^2(M)$ eigenfunction expansion given by
\begin{equation*}%\label{e17}
	p_M(x,y;t)=\sum_{j=1}^{\infty}e^{-t\mu_j(M)}u_{j,M}(x)u_{j,M}(y).
\end{equation*}
It follows by Fubini's theorem that
\begin{equation}\label{e18}
	H_{\Om}(t)=\int_{\Om}dx\,\int_{\Om}dy\, p_M(x,y;t)=\sum_{j=1}^{\infty}e^{-t\mu_j(M)}\Big(\int_{\Om}u_{j,M}\Big)^2.
\end{equation}
We have by \eqref{e18} that
\begin{equation*}%\label{e19}
	\lim_{t\rightarrow\infty}H_{\Om}(t)=\Big(\int_{\Om}u_{1,M}\Big)^2=\frac{|\Om|^2}{|M|}<|\Om|,
\end{equation*}
by hypothesis. Since $H_{\Om}(0)=|\Om|$ we conclude that $H_{\Om}(t)$ is not constant. Then
\begin{align*}%\label{e20}
	\frac{d}{dt}H_{\Om}(t)&=-\sum_{j=1}^{\infty}\mu_j(M)e^{-t\mu_j(M)}\Big(\int_{\Om}u_{j,M}\Big)^2\nonumber\\&
	=-\sum_{j=2}^{\infty}\mu_j(M)e^{-t\mu_j(M)}\Big(\int_{\Om}u_{j,M}\Big)^2,
\end{align*}
is not identically equal to the $0$-function. Let
\begin{equation*}%\label{e21}
	j^*=\min\{j\in\N:j\ge 2,\,\int_{\Om}u_{j,M}\ne 0 \}.
\end{equation*}
Then
\begin{equation*}%\label{e22}
	\frac{d}{dt}H_{\Om}(t)\le -\mu_{j^*}(M)e^{-t\mu_j^*(M)}\Big(\int_{\Om}u_{j^*,M}\Big)^2<0,
\end{equation*}
and $t\mapsto H_{\Om}(t)$ is strictly decreasing.
Similarly
\begin{equation*}%\label{e23}
	\frac{d^2}{dt^2}H_{\Om}(t)\ge \big(\mu_{j^*}(M)\big)^2e^{-t\mu_j^*(M)}\Big(\int_{\Om}u_{j^*,M}\Big)^2>0,
\end{equation*}
and $t\mapsto H_{\Om}(t)$ is strictly convex.

\medskip

\noindent(iv)By \eqref{e4} and \eqref{e5},
\begin{equation*}
	H_{\Om,\psi}(t) =\int_{\Om}\dd y\, u_{\Om,1}(y;t)\psi(y).
\end{equation*}
Since $t\mapsto u_{\Om,1}(y;t)$ is (strictly) decreasing, the integrand in the right-hand side is (strictly) decreasing.
\end{proof}

\section{Proof of Theorem \ref{the3} \label{sec4}}
\begin{proof}
	Since the initial datum is symmetric with respect to the hyperplane $x_1=0$, there is no heat flow across this hyperplane, and the heat equation satisfies Neumann boundary conditions at $x_1=0$.
	The Neumann heat kernel for the half-space $\{x\in\R^m:x_1>0\}=\R_{+}^m$ is denoted and given by
	\begin{equation}\label{ba3}
		\pi_{\R_{+}^m}(x,y;t)=(4\pi t)^{-m/2}\Big(e^{-|x-y|^2/(4t)}+e^{-|x+y|^2/(4t)}\Big).
	\end{equation}
	Hence the solution of \eqref{e1} with $M=\Om_{\delta}$ and $\psi={\bf 1}_{\Om_{\delta}}$ for $x_1>0$ is given by
	\begin{equation}\label{ba4}
		u_{\Om_{\delta}}(x;t)=\int_{B_{\delta}(c_2)} \dd y\, \pi_{\R_{+}^m}(x,y;t),\,x_1>0.
	\end{equation}
	To show that $\Omega_{\delta}$ is a strictly decreasing temperature set we have to show that
	\begin{equation*}%\label{ba5}
		\frac{\partial u_{\Om_{\delta}}(x;t)}{\partial t}<0,\, t>0,\,x\in B_{\delta}(c),
	\end{equation*}
	where we have put $c=c_2$. By \eqref{ba3} and \eqref{ba4} we find that
	\begin{align}\label{ba6}
		\frac{\partial u_{\Om_{\delta}}(x;t)}{\partial t}&=\frac{1}{2t(4\pi t)^{m/2}}\int_{B_{\delta}(c)}\dd y\,\nonumber\\&\Big(e^{-|x-y|^2/(4t)}\Big(\frac{|x-y|^2}{2t}-m\Big)+e^{-|x+y|^2/(4t)}\Big(\frac{|x+y|^2}{2t}-m\Big)\Big).
	\end{align}
	We show that the integrand in the right-hand side of \eqref{ba6} is strictly negative for all
	\begin{equation}\label{ba7}
		x\in B_{\delta}(c),\,y\in B_{\delta}(c),\,t\ge 4\delta^2.
	\end{equation}
	This in turn implies that the left-hand side of \eqref{ba6} is strictly negative for all $t\ge 4\delta^2$.
	By \eqref{ba7},
	\begin{equation}\label{ba77}
		|x-y|\le 2\delta,\, 2\le |x+y|\le 2(1+2\delta).
	\end{equation}
	Hence
	\begin{equation}\label{ba8}
		e^{-\delta^2/t}\le e^{-|x-y|^2/(4t)}\le 1,\,e^{-(1+2\delta)^2/t}\le e^{-|x+y|^2/(4t)}\le e^{-1/t}.
	\end{equation}
	Hence by \eqref{ba7}, \eqref{ba8}, \eqref{ba77} and $\delta=\frac{1}{20}$ we obtain for $t\ge 4\delta^2$,
	\begin{align}\label{ba9}
		e^{-|x-y|^2/(4t)}&\Big(\frac{|x-y|^2}{2t}-m\Big)+e^{-|x+y|^2/(4t)}\Big(\frac{|x+y|^2}{2t}-m\Big)\nonumber\\&
		\le \frac{2\delta^2}{t}-me^{-\delta^2/t}+\frac{2(1+2\delta)^2}{t}e^{-1/t}-me^{-(1+2\delta)^2/t}\nonumber\\&
		\le \frac{2\delta^2}{t}-e^{-\delta^2/t}+\frac{2(1+2\delta)^2}{t}e^{-1/t}-e^{-(1+2\delta)^2/t}\nonumber\\&
		\le \frac{0.005}{t}-e^{-1/4}+\frac{2.42}{t}e^{-1/t}-e^{-1.21/t}.
	\end{align}
	
	\smallskip
	First we consider the case $t\ge 4\delta^2=0.01$. Since $t\mapsto t^{-1}e^{-1/t}$ is strictly decreasing for $t>1$ and since $t\mapsto \frac{0.005}{t}$ and $t\mapsto -e^{-1.21/t}$ are strictly decreasing for all $t>0$ it remains to show that the right-hand side of \eqref{ba9} is strictly less than $0$ on the interval $[0.01,1]$. To reduce this interval further we consider the case $0.01\le t\le \frac14$. On this interval we have that the first term in the right-hand side of \eqref{ba9} is bounded from above by $\frac12$, and that the third term in the right-hand side is bounded from above by $9.68e^{-4}$. One verifies, by for example using Wolfram Alpha \cite{WA}, that $\frac12+ 9.68e^{-4}-e^{-1/4}<-0.1$. Next we consider the interval $[\frac14,1]$.
	On that interval the first term in the right-hand side of \eqref{ba9} is bounded from above by $\frac{1}{50}$. We have that
	\begin{align}\label{ba78}
		\max_{\frac14\le t\le 1}\big(\frac{2.42}{t}e^{-1/t}&-e^{-1.21/t}\big)\nonumber\\&\le \max_{\frac14\le t\le 1}\big(\frac{2.42}{t}-1\big)e^{-1/t}+\max_{\frac14\le t\le 1}\Big(e^{-1/t}-e^{-1.21/t}\Big)\nonumber\\&
		\le \max_{t\ge 0}\big(\frac{2.42}{t}-1\big)e^{-1/t}+\max_{0<t\le 1}\Big(e^{-1/t}-e^{-1.21/t}\Big).
	\end{align}
	It is elementary to verify that the maximum in the first term of the right-hand side of \eqref{ba78} is attained at $t=\frac{121}{171}$. This gives
	\begin{equation*}%\label{ba79}
		\max_{t\ge 0}\big(\frac{2.42}{t}-1\big)e^{-1/t}=\frac{121}{50}e^{-171/121}.
	\end{equation*}
	Furthermore $t\mapsto e^{-1/t}-e^{-1.21/t}$ is increasing on $[0,1]$. Hence
	\begin{equation*}%\label{ba80}
		\max_{0<t\le 1}\Big(e^{-1/t}-e^{-1.21/t}\Big)=e^{-1}-e^{-1.21}.
	\end{equation*}
	One verifies, by for example using Wolfram Alpha \cite{WA}, that
	\begin{equation*}%\label{ba81}
		\frac{1}{50}+\frac{121}{50}e^{-171/121}+e^{-1}-e^{-1.21}-e^{-1/4}<-0.10019.
	\end{equation*}

	\smallskip

	Next we consider the case $0<t\le 4\delta^2,\, \delta=\frac{1}{20}$. For the second term in the right-hand side of \eqref{ba6} we have by \eqref{ba8} and \eqref{ba77} that
	\begin{equation*}%\label{ba10}
		e^{-|x+y|^2/(4t)}\Big(\frac{|x+y|^2}{2t}-m\Big)\le \frac{2(1+2\delta)^2}{t}e^{-1/t}.
	\end{equation*}
	Hence
	\begin{equation*}%\label{ba11}
		\int_{B_{\delta}(c)}dy\,e^{-|x+y|^2/(4t)}\Big(\frac{|x+y|^2}{2t}-m\Big)\le \omega_m\delta^m\frac{2(1+2\delta)^2}{t}e^{-1/t},
	\end{equation*}
	and the second term in the right-hand side of \eqref{ba6} is bounded from above by
	\begin{equation}\label{ba12}
		\frac{\omega_m\delta^m}{(4\pi t)^{m/2}}\frac{(1+2\delta)^2}{t^2}e^{-1/t}.
	\end{equation}
	%to $\frac{\partial u_{\Om_{\delta}}(x;t)}{\partial t}$.
	
	To bound the first term in the right-hand side of \eqref{ba6} we rewrite this term as
	\begin{align*}%\label{ba13}
		\frac{1}{2t(4\pi t)^{m/2}}\int_{B_{\delta}(c)}\dd y\Big(e^{-|x-y|^2/(4t)}\Big(\frac{|x-y|^2}{2t}-m\Big)=\frac{\partial}{\partial t}\int_{B_{\delta}(c)}\dd y\,p_{\R^m}(x,y;t).
	\end{align*}
	Estimating this term for all $x\in B_{\delta}(c)$ is equivalent to estimating
	\begin{equation}\label{ba14}
		\frac{\partial}{\partial t}\int_{B_{\delta}(0)}\dd y\,p_{\R^m}(x,y;t),\, x\in B_{\delta}(0).
	\end{equation}
	Since $B_{\delta}(0)$ is convex we have by \cite[Proposition 8]{MvdB1} that the expression under \eqref{ba14} is strictly negative. Below we quantify this derivative as follows.
	Changing the variable $y-x=t^{1/2}\eta$ yields
	\begin{align*}%\label{ba15}
		\frac{\partial}{\partial t}\int_{B_{\delta}(0)}\dd y\,p_{\R^m}(x,y;t)&=\frac{1}{(4\pi)^{m/2}}\frac{\partial}{\partial t}\int_{B_{\delta t^{-1/2}}(-x)}\dd \eta\,e^{-\eta^2/4}\nonumber\\&
		=-\frac{\delta}{2(4\pi)^{m/2}t^{3/2}}\frac{\partial}{\partial \rho}\int_{B_{\rho}(-x)}\dd \eta\,e^{-\eta^2/4}\Big|_{\rho=\delta t^{-1/2}}\nonumber\\&
		=-\frac{\delta}{2(4\pi)^{m/2}t^{3/2}}\frac{\partial}{\partial \rho}\int_{B_{\rho}(0)}\dd \eta\,e^{-|\eta-x|^2/4}\Big|_{\rho=\delta t^{-1/2}}\nonumber\\&
		=-\frac{\delta}{2(4\pi)^{m/2}t^{3/2}}\int_{\partial B_{\rho}(0)}\mathcal{H}^{m-1}(\dd \eta)\,e^{-|\eta-x|^2/4}\Big|_{\rho=\delta t^{-1/2}},
	\end{align*}
	where $\mathcal{H}^{m-1}(\dd \eta)$ denotes the surface measure. For $x\in B_{\delta t^{-1/2}}(0)$ and $\eta\in\partial B_{\delta t^{-1/2}}(0)$ we have that
	$|x-\eta|^2\le 4\delta^2/t$. This gives
	\begin{align}\label{ba16}
		\frac{\partial}{\partial t}\int_{B_{\delta}(0)}\dd y\,p_{\R^m}(x,y;t)&\le -\frac{\delta}{2(4\pi)^{m/2}t^{3/2}}\int_{\partial B_{\delta t^{-1/2}}(0)}\mathcal{H}^{m-1}(\dd \eta)\,e^{-\delta^2/t}\nonumber\\&
		=-\frac{m\omega_m\delta^m}{2t(4\pi t)^{m/2}}e^{-\delta^2/t}.
	\end{align}
	By \eqref{ba6}, \eqref{ba12} and \eqref{ba16},
	\begin{equation}\label{ba17}
		\frac{\partial u_{\Om_{\delta}}(x;t)}{\partial t}\le -\frac{m\omega_m\delta^m}{2t(4\pi t)^{m/2}}e^{-\delta^2/t}+\frac{\omega_m\delta^m}{(4\pi t)^{m/2}}\frac{(1+2\delta)^2}{t^2}e^{-1/t},\,x\in B_{\delta}(0).
	\end{equation}
	For $0<t\le 4\delta^2,\,\delta=\frac{1}{20}$, since $t \mapsto te^{(1-\delta^2)t^{-1}}$ is decreasing, we have that
\begin{equation*}%\label{ba19}
		te^{(1-\delta^2)t^{-1}}\ge 4\delta^2e^{(1-\delta^2)/(4\delta^2)}>2(1+2\delta)^2\ge\frac{2}{m}(1+2\delta)^2.
	\end{equation*}
This implies that the right-hand side of \eqref{ba17} is strictly negative.
\end{proof}

\section{Proof of Theorem \ref{the5} \label{sec6}}
\begin{proof}
	A straightforward calculation shows by \eqref{ee} that
	\begin{equation}\label{e53}
		\frac{\partial^2p_{\R^m}(x,y;t)}{\partial t^2}= \frac{1}{t^2}\Big(\Big(\frac{m+2}{2}-\frac{b}{t}\Big)^2-\frac{m+2}{2}\Big)p_{\R^m}(x,y;t),
	\end{equation}
	where
	\begin{equation*}%\label{e54}
		b=\frac14|x-y|^2\le \frac14(\textup{diam}(\Om))^2.
	\end{equation*}
	For all $t\ge (\textup{diam}(\Om))^2$ we have $\frac{b}{t}\le \frac14$. This, together with \eqref{e53} gives
	\begin{equation*}%\label{e55}
		\frac{\partial^2p_{\R^m}(x,y;t)}{\partial t^2}\ge \frac{1}{t^2}\Big(\frac{4m^2+4m-7}{16}\Big)p_{\R^m}(x,y;t),\, t\ge (\textup{diam}(\Om))^2.
	\end{equation*}
	Integrating both sides with respect to $x\in\Om,y\in\Om$ gives the assertion under (i).
	
	\smallskip
	
	To prove (ii) we note that, uniformly in $x$ and $y$ in $\Om$,
	\begin{equation*}%\label{e56}
		p_{\R^m}(x,y;t)\ge \frac{e^{-1/4}}{(4\pi t)^{m/2}},\,t\ge (\textup{diam}(\Om))^2.
	\end{equation*}
	Integrating both sides with respect to $x\in\Om,y\in\Om$ gives
	\begin{equation*}%\label{e57}
		H_{\Om}(t)\ge e^{-1/4}\frac{|\Om|^2}{(4\pi t)^{m/2}}.
	\end{equation*}
	This, together with \eqref{e49}, gives the assertion under (ii).
	
	\smallskip
	
	To prove (iii) we first consider $t\ge (\textup{diam}(\Om))^2$, and integrate \eqref{e50} between $s$ and $\infty$ where $s \geq (\textup{diam}(\Om))^2$.
	This gives the second inequality in \eqref{e51}.
	Since the heat content is convex, its first derivative is increasing and continuous. This proves the first inequality in \eqref{e51}.
	
	\smallskip
	
	To prove (iv) we have that
	\begin{equation*}%\label{e58}
		p_{\R^m}(x,y;t)\le (4\pi t)^{-m/2}
	\end{equation*}
	implies
	\begin{equation*}%\label{e59}
		H_{\Om}(t)\le \frac{|\Om|^2}{(4\pi t)^{m/2}}.
	\end{equation*}
	This, together with (iii), yields (iv).
\end{proof}

\section{Proof of Theorem \ref{the2} \label{sec3}}
\begin{proof}
We define the heat loss of $\Omega_c$ at $t$ by
\begin{equation}\label{e27}
F_{\Omega_c,\psi}(t)=H_{\Omega_c,\psi}(0)-H_{\Omega_c,\psi}(t).
\end{equation}
It follows that
\begin{align*}%\label{e28}
F_{\Omega_c,\psi}(t)&=\int_{\R^m}\dd x\,\int_{\Omega_c}\dd y\,\psi(y)p_{\R^m}(x,y;t)-\int_{\Omega_c}\dd x\,\int_{\Omega_c}\dd y\,\psi(y)p_{\R^m}(x,y;t)\nonumber\\&
=\int_{\R^m\setminus\Omega_c}\dd x\,\int_{\Omega_c}\dd y\,\psi(y)p_{\R^m}(x,y;t).
\end{align*}
We will show that if $c$ satisfies \eqref{e25} then $F_{\Om_c,\psi}(1)>F_{\Om_c,\psi}(t^*)>0$. Since $t^*>1$ and $F_{\Omega_c,\psi}(0)=0$ we infer that the heat loss, and hence the heat content, is not monotone.

We consider $\mathbb{R}^m \setminus \Omega_c$ and, for $r_2 > r_1 > 0$, define $A_{[r_1,r_2]}=\{x\in\R^m:r_1\le|x|\le r_2\}$
	and $A_{[r_1, \infty)} =  \{x\in\R^m: |x| \ge r_1\}$.
For $x\in  A_{[1,2]}$ and $y\in B_1{ (0)}$, $|x-y|^2<9$. This gives
\begin{align*}%\label{e30}
F_{\Om_c,\psi}(t)&\ge \int_{A_{[1,2]}}\dd x\,\int_{B_1{ (0)}} \dd y\,(4\pi t)^{-m/2}e^{-9/(4t)}\nonumber\\&
=(4\pi t)^{-m/2}e^{-9/(4t)}\omega_m^2(2^m-1).
\end{align*}
Hence
\begin{equation}\label{e31}
F_{\Omega_c,\psi}(1)\ge e^{-9/4}\frac{\omega_m^2(2^m-1)}{(4\pi)^{m/2}}.
\end{equation}
To obtain an upper bound for $F_{\Omega_c,\psi}$ we have
\begin{align}\label{e32}
\int_{\R^m\setminus\Omega_c}\dd x\,\int_{A_{(2,c)}}\dd y\,\psi(y)p_{\R^m}(x,y;t)=0,
\end{align}
and
\begin{align}\label{e33}
\int_{A_{[1,2]}}\dd x\,\int_{B_1{(0)}}\dd y\,p_{\R^m}(x,y;t)&\le (4\pi t)^{-m/2}|B_1{ (0)}||A_{[1,2]}|\nonumber\\&
=(4\pi t)^{-m/2}\omega_m^2(2^m-1).
\end{align}
For $|x|\ge c$ and $y\in B_1{ (0)}$, we have that $|y|<|x|/2$. Hence
\begin{align}\label{e35}
\int_{A_{[c,\infty)}}&\dd x\,\int_{B_1{ (0)}}\dd y\,p_{\R^m}(x,y;t)\le (4\pi t)^{-m/2}\int_{A_{[c,\infty)}}dx\,\int_{B_1 (0)}dy\,e^{-|x|^2/(16t)}\nonumber\\&
\le (4\pi t)^{-m/2}e^{-c^2/(32t)}\int_{\R^m}\dd x\,\int_{B_1{ (0)}}\dd y\,e^{-|x|^2/(32t)}\nonumber\\&
\le 2^{3m/2}\omega_me^{-c^2/(32t)}.
\end{align}
Putting \eqref{e33} and \eqref{e35} together yields
\begin{align}\label{e37}
\int_{\R^m\setminus\Omega_c}\dd x\,&\int_{B_1{ (0)}}\dd y\,p_{\R^m}(x,y;t)\nonumber\\&\le\frac{\omega_m^2(2^m-1)}{(4\pi t)^{m/2}}+2^{3m/2}\omega_me^{-c^2/(32t))}.
\end{align}
Combining \eqref{e32} and \eqref{e37} gives
\begin{equation}\label{e38}
\int_{\R^m\setminus\Omega_c}\dd x\,\int_{\Omega_c}\dd y\,\psi(y)p_{\R^m}(x,y;t)\le \frac{\omega_m^2(2^m-1)}{(4\pi t)^{m/2}}+2^{3m/2}\omega_me^{-c^2/(32t)}.
\end{equation}
If, for $t^* > e^{9/(2m)}$, \eqref{e25} holds, then
\begin{equation*}%\label{e40}
	2^{3m/2}e^{-c^2/(32t^*)}<\frac{\omega_m(2^m-1)}{(4\pi)^{m/2}}\big(e^{-9/4}-{t^*}^{-m/2}\big),
\end{equation*}
which implies that the right-hand side of \eqref{e38} is bounded from above by the right-hand side of \eqref{e31} as required.
\end{proof}

\section{Proof of Theorem \ref{the4} \label{sec5}}

\begin{proof}
The set $B_1(0)$ is convex and $\psi\ge 0$ and measurable. It follows by Theorem \ref{the1}(iv) that the heat content is strictly decreasing.  Moreover, the heat content is strictly positive and decreasing to $0$.
To prove non-convexity it therefore suffices to show that the right-derivative $H'^+_{B_1(0),\psi}(0)=0$. We have the following
\begin{align*}%\label{e44}
-H'^+_{B_1(0),\psi}(0)&=\lim_{t\downarrow 0}t^{-1}\big(H_{B_1(0),\psi}(0)- H_{B_1(0),\psi}(t)\big)\nonumber\\&
=\lim_{t\downarrow 0}t^{-1}F_{B_1(0),\psi}(t)\nonumber\\&
=\lim_{t\downarrow 0}\frac{1}{t}\int_{\{|x|>1\}}\dd x\,\int_{B_1(0)}\dd y\,(4\pi t)^{-m/2}e^{-\frac{|x-y|^2}{4t}}(1-|y|)^{\alpha},
\end{align*}
where we have used \eqref{e27}. By the radial symmetry of $\psi$, the map
\begin{equation*}%\label{e45}
x\mapsto \int_{B_1(0)}\dd y\,(4\pi t)^{-m/2}e^{-\frac{|x-y|^2}{4t}}(1-|y|)^{\alpha}
\end{equation*}
is radially symmetric, and depends only on $|x|$. Without loss of generality we put $x=vr$, where $v=(1,0, \dots ,0).$  Changing variable  $y=v-\eta$ yields
\begin{align}\label{e46}
&-H'^+_{B_1(0),\psi}(0)\nonumber\\&=\lim_{t\downarrow 0}\frac{m\omega_m}{t}\int_{(1,\infty)}\dd r\,r^{m-1}\int_{\{|v-\eta|<1\}}\dd \eta\,(4\pi t)^{-m/2}e^{-\frac{|v(r-1)+\eta|^2}{4t}}(1-|v-\eta|)^{\alpha}.
\end{align}
Note that $\eta$ is a vector with $|v\cdot\eta|>0$. Hence $e^{-\frac{|v(r-1)+\eta|^2}{4t}}\le e^{-\frac{(r-1)^2+|\eta|^2}{4t}}$.
Furthermore, $|v-\eta|\le |v|+|\eta|=1+|\eta|$. So $1-|v-\eta|\ge -|\eta|$. Also $|v-\eta|\ge |v|-|\eta|=1-|\eta|$. Hence $1-|v-\eta|\le |\eta|$.
Hence the right-hand side of \eqref{e46} is bounded from above by
\begin{align*}%\label{e47}
&\lim_{t\downarrow 0}\frac{m\omega_m}{t}\int_{(1,\infty)}\dd r\,r^{m-1}\int_{\{|v-\eta|<1\}}\dd \eta\,(4\pi t)^{-m/2}e^{-\frac{(r-1)^2+|\eta|^2}{4t}}|\eta|^{\alpha}\nonumber\\&\le
\lim_{t\downarrow 0}\frac{m\omega_m}{t}\int_{(1,\infty)}\dd r\,r^{m-1}(4\pi t)^{-m/2}e^{-\frac{(r-1)^2}{4t}}\int_{\R^m}\dd \eta\,|\eta|^{\alpha}e^{-|\eta|^2/(4t)}\nonumber\\&
=\lim_{t\downarrow 0}\frac{(m\omega_m)^2}{2t}(4t)^{(m+\alpha)/2}\Gamma((m+\alpha)/2)(4\pi t)^{-m/2}\int_{(0,\infty)}\dd r\,(1+r)^{m-1}e^{-\frac{ r^2}{4t}}\nonumber\\&
\leq \lim_{t\downarrow 0}\frac{(m\omega_m)^2}{2t}(4t)^{(m+\alpha)/2}\Gamma((m+\alpha)/2)(4\pi t)^{-m/2}\nonumber\\&
\hspace{45mm}\times\int_{(0,\infty)}\dd r\,2^{m-1}(1+r^{m-1})e^{-\frac{ r^2}{4t}}\nonumber \\&
=\lim_{t\downarrow 0}\big(C_{m,\alpha,1}t^{(\alpha-1)/2}+C_{m,\alpha,2}t^{(\alpha-2+m)/2}\big)=0,
\end{align*}
since $\alpha>1$, and $C_{m,\alpha,1}<\infty$ and $C_{m,\alpha,2}<\infty$ are  strictly finite constants depending on $m$ and on $\alpha$ only.
\end{proof}

\noindent {\bf Acknowledgements.}\, It is a pleasure to thank Dorin Bucur for some helpful discussions.
We are grateful to the referees for their helpful comments.

\bigskip

\noindent  {\bf ORCID.}  Michiel van den Berg: http://orcid.org/0000-0001-7322-6922\\

\noindent  {\bf ORCID.}  Katie Gittins: https://orcid.org/0000-0001-8136-8437\\

\end{document}